\documentclass[12pt]{amsart}
\usepackage{amssymb,amsfonts,amsmath,amsopn,
amstext,amscd,latexsym,amsthm,enumerate,mathrsfs,bbm,color}

\usepackage[margin=36mm]{geometry}
\usepackage[all,cmtip]{xy}
\usepackage{longtable}

\newtheorem{theorem}{Theorem}
\newtheorem{thm}[theorem]{Theorem}

\theoremstyle{remark}

\numberwithin{equation}{section}

\begin{document}

\title[Monomial, monolithic characters]
{Solvable groups whose monomial, monolithic characters have prime power codegrees}

\author[X. Chen] {Xiaoyou Chen}
\address{School of Sciences, Henan University of Technology, Zhengzhou 450001, China}
\email{cxymathematics@hotmail.com}

\author[M. L. Lewis]{Mark L. Lewis}
\address{Department of Mathematical Sciences, Kent State University, Kent, OH 44242, USA}
\email{lewis@math.kent.edu}


\date{\today}


\begin{abstract}
In this note,
we prove that if $G$ is solvable and ${\rm cod}(\chi)$ is a $p$-power for every nonlinear, monomial, monolithic $\chi\in {\rm Irr}(G)$ or every nonlinear, monomial, monolithic $\chi \in {\rm IBr} (G)$,
then $P$ is normal in $G$, where $p$ is a prime and $P$ is a Sylow $p$-subgroup of $G$.
\end{abstract}

\maketitle


All groups in this paper are finite.  We refer the reader for notation to \cite{Isaacs1}.  Let $G$ be a group, and we write ${\rm Irr} (G)$ for the set of irreducible characters of $G$.  For a fixed prime $p$, we use ${\rm IBr} (G)$ to denote the set of irreducible $p$-Brauer characters of $G$.

A number of papers have recently considered the following general question: 
let $G$ be a solvable group and let $\mathcal {A}$ be a fixed subset of either ${\rm Irr} (G)$ or ${\rm IBr} (G)$; if all the elements in $\mathcal {A}$ have a given property, what can be said about the structure of $G$?  In \cite{lu}, \cite{PaLu}, and \cite{PaLu2}, the set ${\mathcal A}$ considered is a subset of nonlinear monomial characters.  Whereas in \cite{lql}, they consider the set $\mathcal {A}$ of nonlinear monolithic characters.  For characters in the given subset, they may further consider those whose degrees are divisible by a prime $p$ or not divisible by a prime $p$.  The properties they look at whether the group include nilpotent, $p$-nilpotent, and $p$-closed.  Recall that a character is {\it monomial} if it induced from a linear character of a subgroup.  A character $\chi$ is said to be  {\it monolithic} if $G/\ker (\chi)$ has a unique minimal normal subgroup.
A group is {\it $p$-closed} if it has a normal Sylow $p$-subgroup and it is {\it $p$-nilpotent} if it has a normal $p$-complement.
In \cite{CL1}, \cite{CL3}, \cite{CL4}, and \cite{CL5}, similar questions are considered for similar subsets of the Brauer characters and the characters in \cite{CY}.

Given a character $\chi \in \mathcal {A}$, we define the {\it codegree} of $\chi$ by $\displaystyle {\rm cod}(\chi) = {|G: \ker\chi|} /{\chi(1)}.$  This usage of the term codegree is due to Qian, Wang and Wei in \cite{Qian1}, although the term codegree of $\chi$ had previously been applied to the quotient $|G|/\chi(1)$ by Chillag, Mann, and Manz in \cite{CMM}.

Let $G$ be a nontrivial group and $p$ be a prime divisor of $|G|$.

Qian for solvable groups in Theorem 1.1 of \cite{Qian2} and Isaacs for nonsolvable groups in \cite{Isaacs2} have proved that for every element $g$ in the group $G$ there exists a character $\chi \in {\rm Irr} (G)$ so that ${\rm cod} (\chi)$ is divisible by every prime divisor of $o(g)$.  It follows that ${\rm cod}(\chi)$ is a $p$-power for every character $\chi \in {\rm Irr}(G)$ if and only if $G$ is a $p$-group.

We now restrict our attention to either nonlinear irreducible characters or nonlinear irreducible Brauer characters that are monomial and monolithic.  We now show that if we assume that the characters or Brauer characters in the given set all have codegrees that are $p$-powers, then $G$ is $p$-closed.  Recall that a group is $p$-closed if it has a normal Sylow $p$-subgroup.

\begin{thm}
Let $G$ be a solvable nontrivial group, and let $p$ be a prime divisor of $|G|$.  Let $\mathcal {A}$ be either the set of nonlinear, monomial, monolithic characters in ${\rm Irr} (G)$ or the set of nonlinear, monomial, monolithic Brauer characters in ${\rm IBr} (G)$.  If ${\rm cod}(\chi)$ is a $p$-power for every $\chi \in \mathcal {A}$, then $G$ is $p$-closed.
\end{thm}

Observe that if $G$ has no nonlinear monomial irreducible characters, then by Theorem 1.1 of \cite{PaLu} that $G$ is abelian.  Also, under the condition that $P$ is a normal Sylow $p$-subgroup in $G$ we do not necessarily have that ${\rm cod} (\chi)$ are $p$-powers for all nonlinear monomial $\varphi\in {\rm Irr}(G)$.  For example, let $G=C_{3}\times Q_{8}$ and $p=3$; then $G$ has a normal Sylow $3$-group, however, ${\rm cod}(\chi)$ are not $3$-powers for the nonlinear, monomial, monolithic characters $\chi \in {\rm Irr}(G)$.

Also, if ${\rm cod}(\chi)$ are $p$-powers for all nonlinear monomial irreducible characters of $G$, we do not have that the Sylow $p$-subgroup $P$ of $G$ is abelian.  For example, let $G={\rm SL}(2, 3)$ and note that ${\rm cod}(\chi)=4$ for the nonlinear monomial character $\chi\in {\rm Irr}({\rm SL}(2,3))$ with degree $3$; however, the Sylow $2$-subgroup of $G$ is not abelian.

Since the following proof is very similar to the proof of Theorem 1 of \cite{CL3}, the solvable part of the proof of Theorem 1.3 of \cite{CL4}, and Theorem 1 of \cite{CL5}, we will suppress many of the details of this proof.  Hence, we will leave out many of the usual details.

\begin{proof}[Proof of Theorem 1]
Suppose that ${\rm cod}(\chi)$ are all $p$-powers for either the characters or the Brauer characters in $\mathcal {A}$.  Let $G$ be a counterexample of minimal order.  

Let $M$ be a minimal normal subgroup of $G$.  Using the inductive hypothesis, we can show that $M$ is a $p'$-group and the unique minimal normal subgroup of $G$.   Applying the Frattini argument, we obtain $G = M {\bf N}_G (P)$ and $M \cap {\bf N}_G (P) = 1$.  Since $M$ is the unique minimal normal subgroup, we obtain ${\bf C}_P (M) = 1$, so $P$ acts faithfully on $M$.  It follows that $P$ acts faithfully on ${\rm Irr}(M)$.  (Since $M$ is a $p'$-group, ${\rm Irr} (M) = {\rm IBr} (M)$.)

It follows by Isaacs' large orbit result (Theorem B of \cite{large}) that there exists a nonprincipal character $\lambda \in {\rm Irr} (M) = {\rm IBr} (M)$ so that $|{\bf C}_{P} (\lambda)| < \sqrt {|P|}$.  Thus, $|P:{\bf C}_{P} (\lambda)| > \sqrt {|P|}$ and then $|P:{\bf C}_{P} (\lambda)|^2 > |P|$.  Write $T$ for the inertia group of $\lambda$ in $G$. Write $S = M {\bf C}_{P} (\lambda)$.  Observe that $S$ is the stabilizer of $\lambda$ in $N=MP$ and $S \leq T$.  In particular, $S  = T \cap N$.  Since $N/M$ is the Sylow $p$-subgroup of $G/M$, it follows that $S/M$ is the Sylow $p$-subgroup of $T/M$. We have that $|N:S| = |P:{\bf C}_P (\lambda)|$ and $|P| = |N:M| = |G:M|_p = |G|_p$, where $|G|_{p}$ denotes the $p$-part of $|G|$.  We deduce that $$|G: T|_{p} = |P:{\bf C}_P (\lambda)| > \sqrt {|P|}.$$

Observe that $M$ is complemented in $T$.  It follows from  a result of Gallagher (see Lemma 1 of \cite{Gallagher}) that there exists some character $\mu\in {\rm Irr} (T)$ such that $\mu_{M} = \lambda$ or a Brauer character $\mu^* \in {\rm IBr} (T)$ so that $\mu^*_M = \lambda$.  By the Clifford correspondence (Theorem 6.11 of \cite{Isaacs1} for characters and Theorem 8.9 of \cite{navbook} for Brauer characters), $\varphi = \mu^{G}\in {\rm Irr}(G)$ and $\varphi^* = (\mu^*)^G$.  Since $\mu$ and $\mu^*$ are linear, we have that $\varphi$ and $\varphi^*$ are monomial and $\varphi (1) = \varphi^* (1) = |G: T|$.

Since $M$ is the unique minimal normal subgroup of $G$, it follows that $\ker\varphi = \ker {\varphi^*}  = 1$.  Otherwise, $M \subseteq \ker\varphi$ and then $M \subseteq \ker\mu$ or $M \subseteq \ker {\varphi^*}$ and $M \subseteq \ker {\mu^*}$, so $\mu_{M}$ or $\mu^*_M$ is the principal character of $M$, which is impossible.  Since $M$ is the unique minimal normal subgroup of $G$, it follows that $\varphi$ and $\varphi^*$ are monolithic.  Notice that $\varphi(1)_{p}^{2} = \varphi^* (1)_p^2 = |G: T|_{p}^{2} > |P| = |G|_{p}$. So $$\varphi(1)_{p}^{2} |G|_{p'} = \varphi^* (1)_{p}^{2} |G|_{p'} >|G|_{p}|G|_{p'}=|G|\geq |G|_{p}.$$

By hypothesis, we have that $${\rm cod} (\varphi) = {\rm cod} (\varphi^*) = \frac{|G: \ker\varphi|} {\varphi(1)} = \frac{|G|} {\varphi(1)}$$
is a $p$-power.  Note that $\varphi(1)^{2} < |G|$ and $(|G|_{p'} \varphi(1)_{p})^{2} < |G|$. Therefore,
$$|G|_{p} > \varphi(1)_{p}^{2} |G|_{p'},$$ a contradiction with the previous paragraph.
\end{proof}


\section*{Acknowledgments}
The first author thanks Cultivation Programme for Young Backbone Teachers in Henan University of Technology,
the Project of Henan Province (182102410049), and the NSFC (11926330, 11926326, 11971189, 11771356).


\end{document}